\documentclass[10pt]{article}
\usepackage{amssymb,amsmath,amsthm,tikz,multirow,nccrules,float,colortbl,arydshln,multicol,ulem,graphicx,subfig}
\usetikzlibrary{arrows,calc}

\title{Cyclotomic points on varieties and all rational $a^3b$-monotiles}
\author{Jinjin Liang, Yixi Liao, Erxiao Wang\thanks{Corresponding author (wang.eric@zjnu.edu.cn).  Research was supported by National Natural Science Foundation of China NSFC-RGC 12361161603 and Key Projects of Zhejiang Natural Science Foundation LZ22A010003.} \\
	Zhejiang Normal University}

\allowdisplaybreaks[4]  

\newcommand\aaa{\alpha}
\newcommand\bbb{\beta}
\newcommand\ccc{\gamma}
\newcommand\ddd{\delta}

\newtheorem{theorem}{Theorem}
\newtheorem{lemma}[theorem]{Lemma}

\newtheorem*{theorem*}{Theorem}

\theoremstyle{definition}
\newtheorem*{definition*}{Definition}
\newtheorem*{case*}{Case}
\newtheorem*{subcase*}{Subcase}

\def\Z{{\mathbb{Z}}}

\numberwithin{equation}{section}

\begin{document}
	\date{}
	
	\maketitle	
\begin{abstract}
	
	By computing all cyclotomic points on some algebraic varieties, we get an independent and efficient way to find all rational $a^3b$-monotiles for the sphere, thereby completing the classification of edge-to-edge monohedral quadrilateral tilings. Both of the previous classifications \cite{lw2} and \cite{cl} depended on many old works of different authors while quite a few  typos and gaps were found.

	{\it2020 MR Subject Classification} Primary 52C20, 05B45; Secondary
	11R18, 11Y50, 14Q25.
		
\end{abstract}

\section{Introduction}
\begin{figure}[htp]
	\centering
	\includegraphics[scale=0.14]{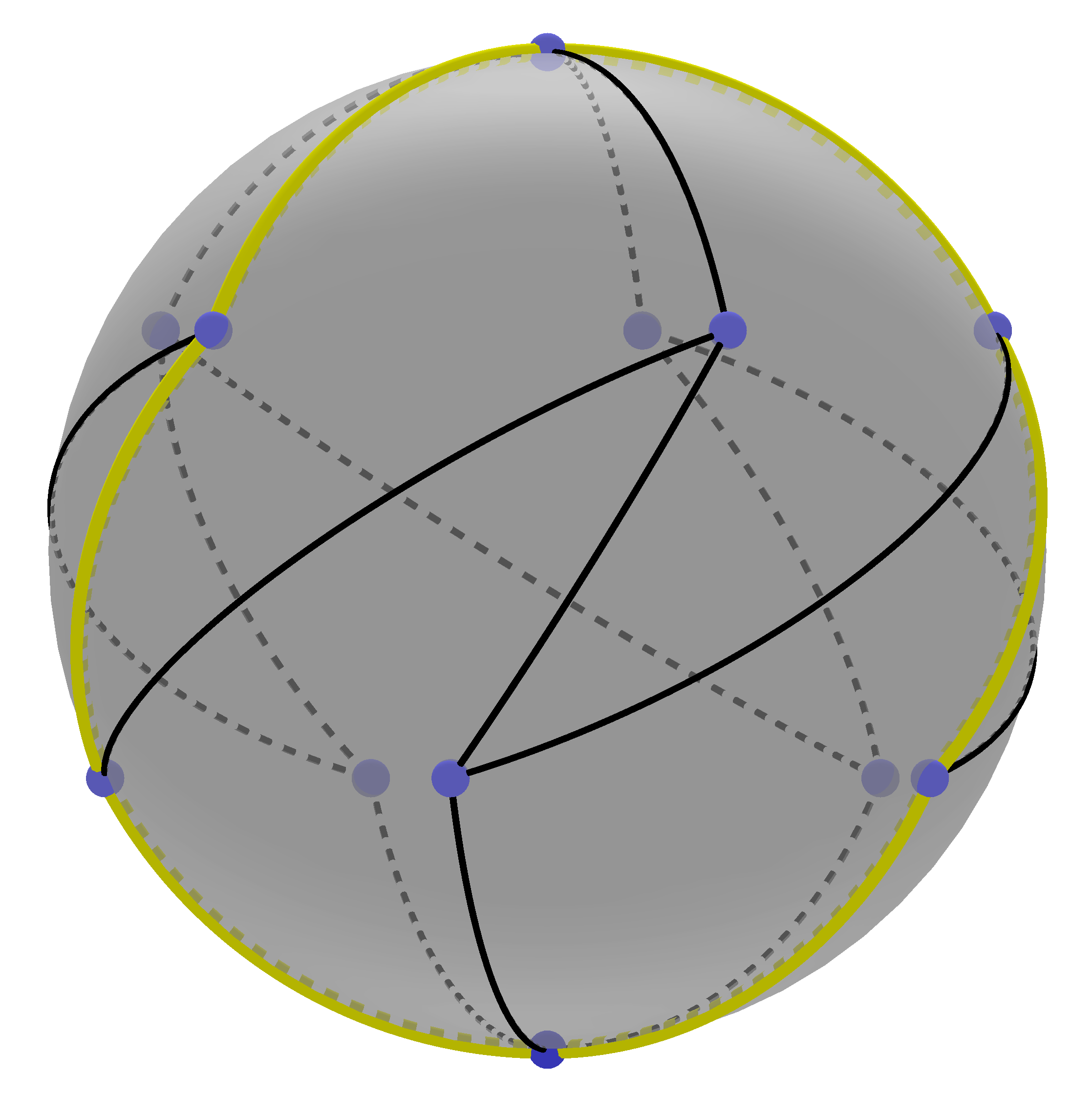}		
	\includegraphics[scale=0.162]{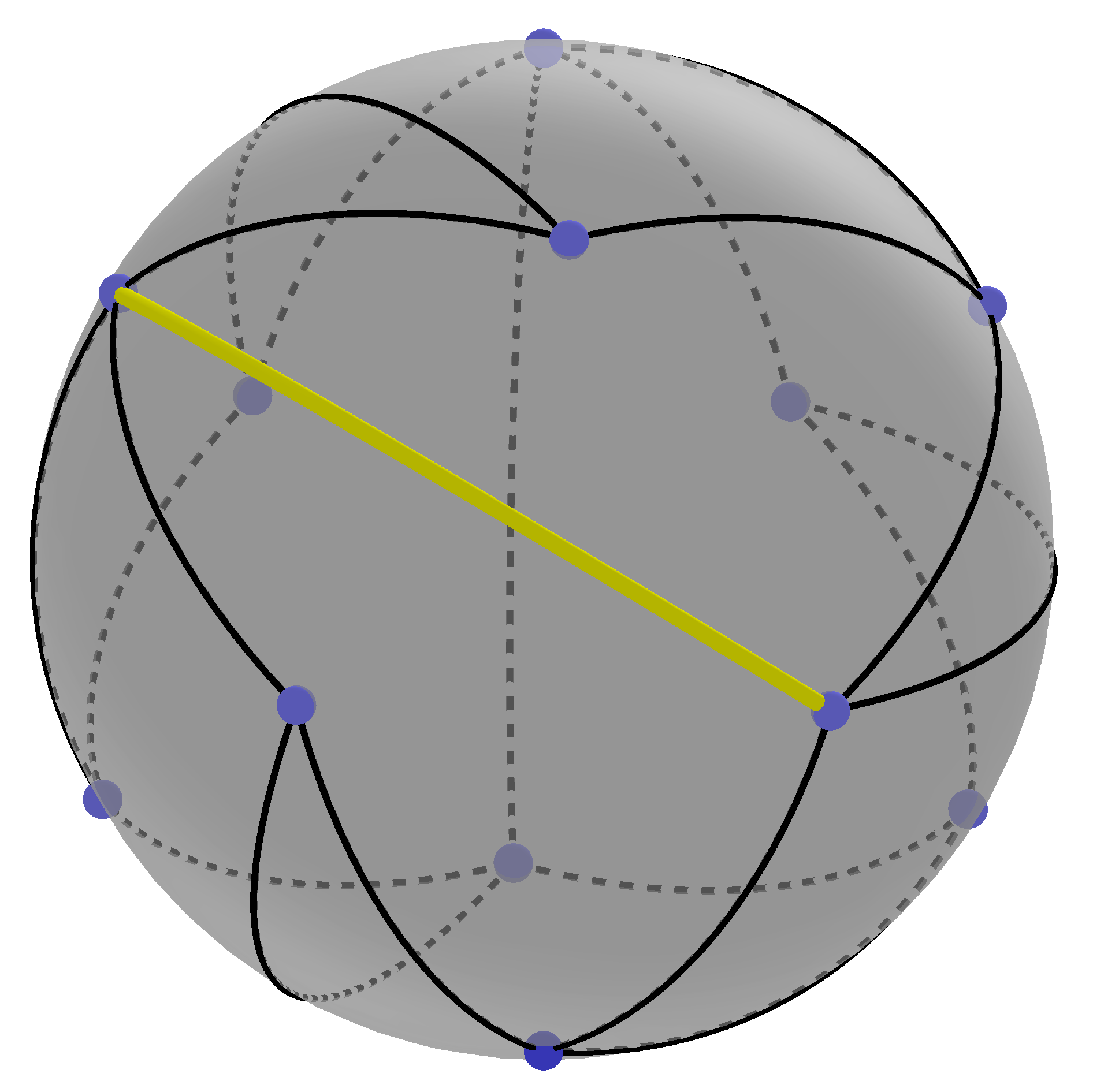} 
	\includegraphics[scale=0.132]{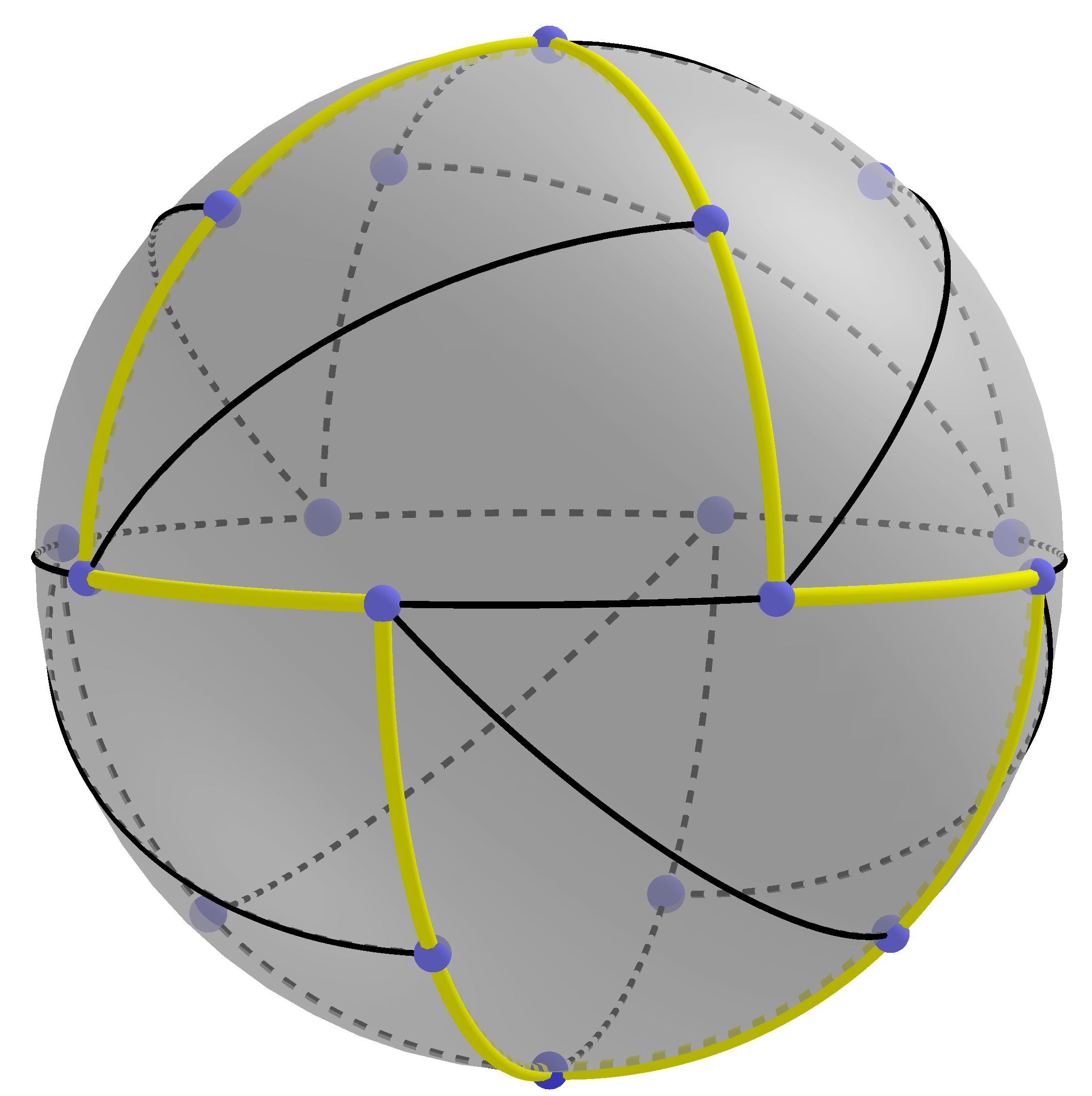}				
	\includegraphics[scale=0.135]{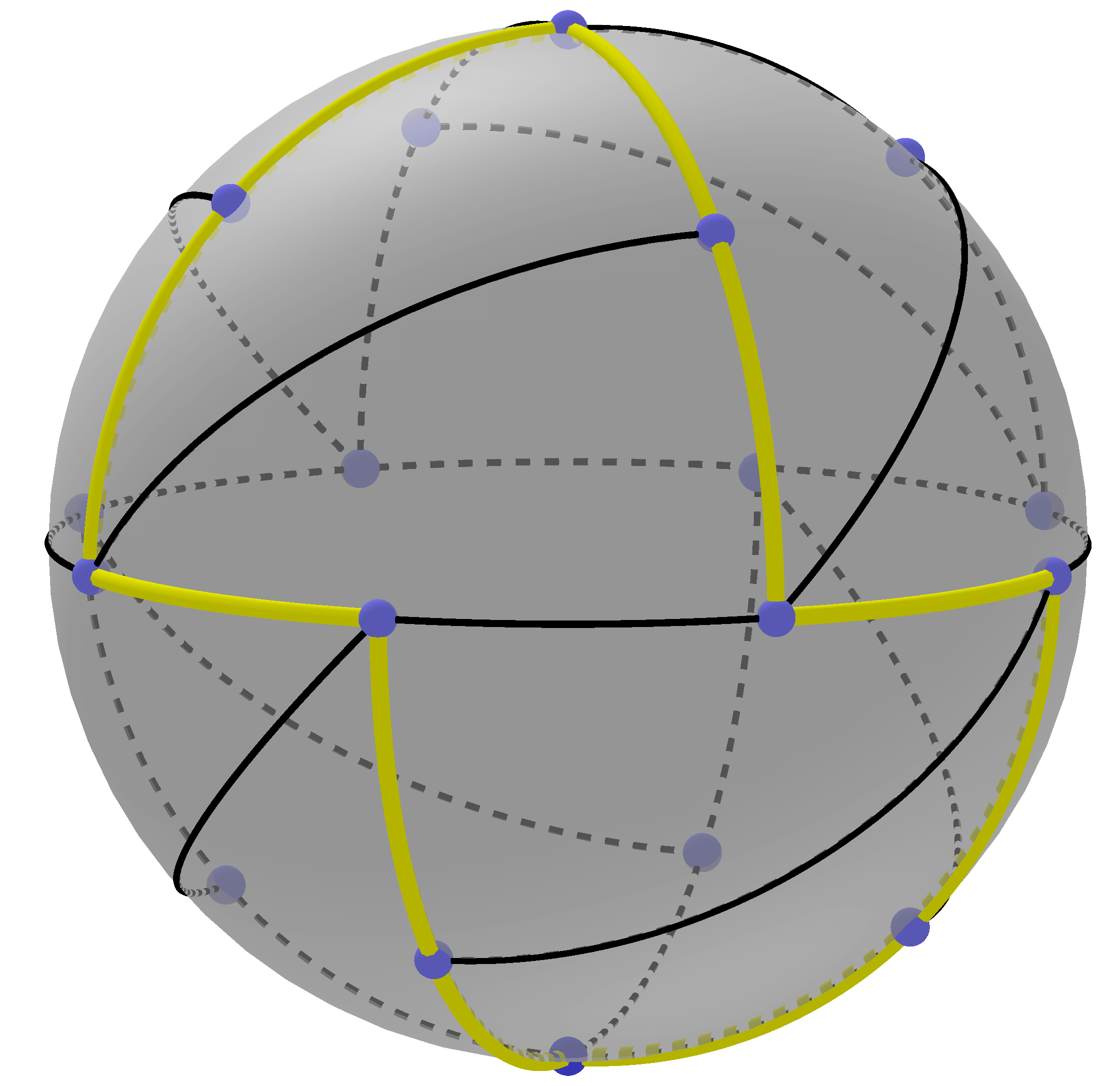}  
	\includegraphics[scale=0.22]{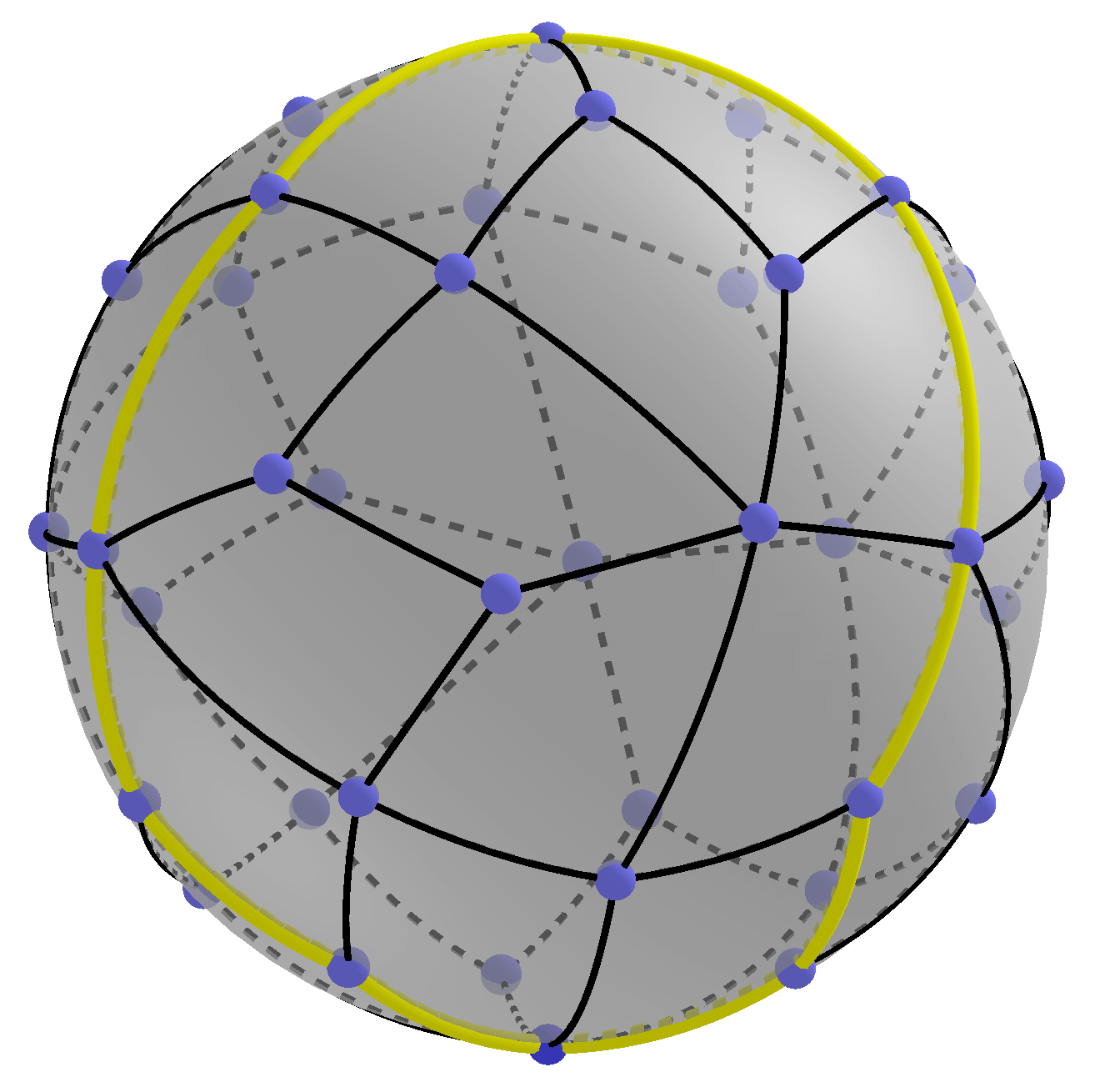}
	\includegraphics[scale=0.20]{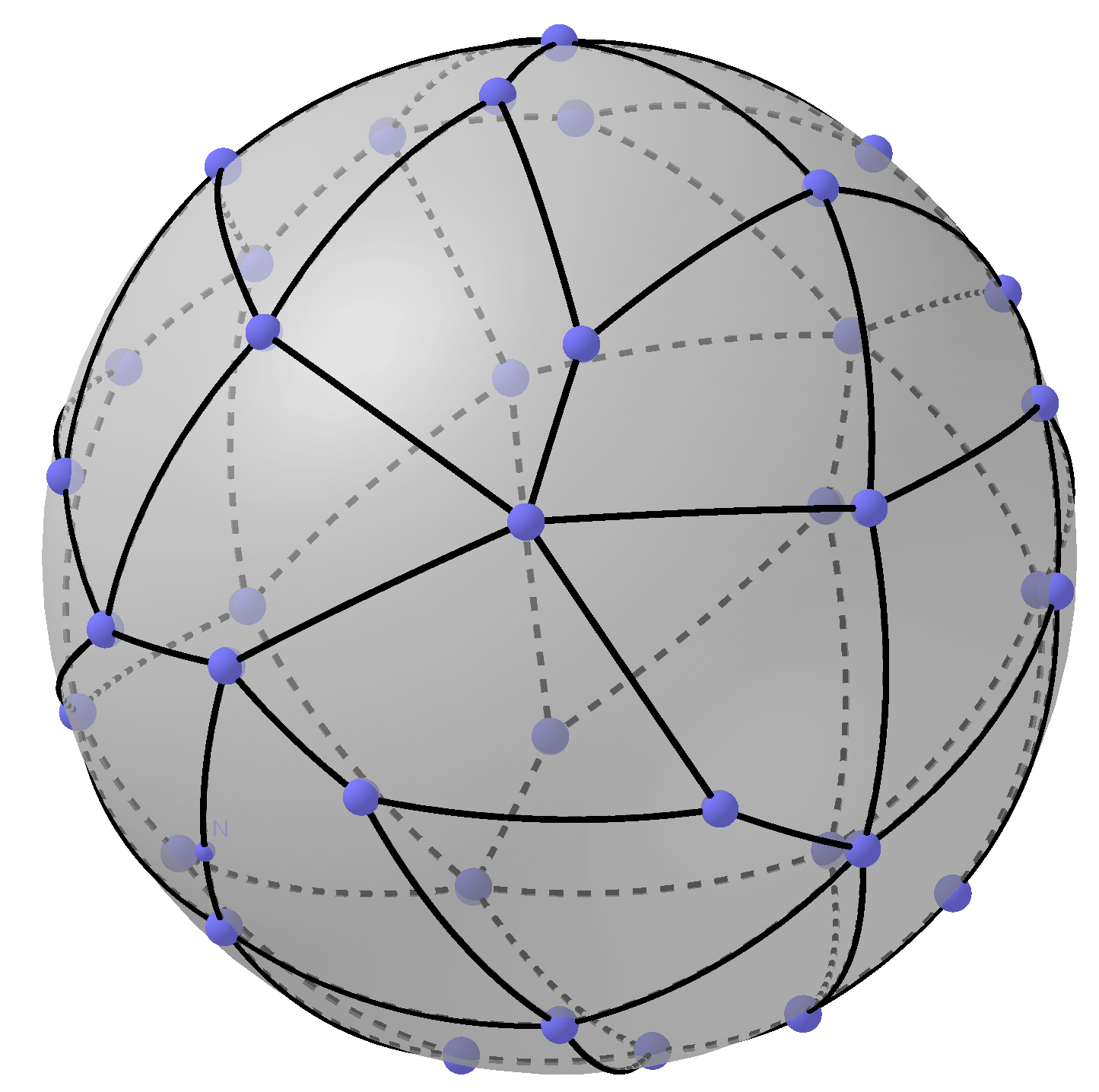}
	\caption{Some ``almost equilateral quadrilateral'' or simply $a^3b$-tilings with $10,10,16,16,36,36$ tiles.} 
	\label{fig 1-1-1}	
\end{figure}

Tilings have numerous applications in real life and science itself. Two great books \cite{ac,gs} introduced their modern  mathematical theory. Open questions like ``which tetrahedra fill space'' \cite{sm} go back to $2300+$ years ago, and led to the first modern study of spherical tilings in 1924 \cite{so}. After 100 years, the classification of edge-to-edge monohedral tilings of the sphere have been completed recently in \cite{wy1,wy2,awy,cly1,lsw,llhw,lw1,lw2,lw3,cly},  mainly by two groups in the Hong Kong University of Science $\&$ Technology and Zhejiang Normal University. Figure \ref{fig 1-1-1} shows six examples of   ``almost equilateral quadrilateral'' or simply $a^3b$-tilings, which together with $a^4b$-pentagonal tilings are the hardest cases in the classification program. 

\begin{figure}[htp]
	\centering
	\begin{tikzpicture}

		\foreach \a in {0}
		{
			\begin{scope}[xshift=3*\a cm]
				\draw
				(-0.8,-0.8) -- (-0.8,0.8) -- (0.8,0.8) -- (0.8,-0.8);

				\draw[line width=1.5]
				(-0.8,-0.8) -- (0.8,-0.8);	    	\end{scope}}   
		
		\draw  (1.5+0.5,0.1)--(2.5+0.5,0.1);

		\draw[line width=1.5] (1.5+0.5,-0.4)--(2.5+0.5,-0.4);
		
		\node at (-0.5,0.5) {\small $\bbb$};
		\node at (0.5,0.5) {\small $\ccc$};
		\node at (-0.5,-0.5) {\small $\aaa$};
		\node at (0.5,-0.5) {\small $\ddd$};
		
		\foreach \a in {0,1,2}
		\node at (90*\a:1) {\small $a$};
		\node at (0,-1.1) {\small $b$};
		
		\node at (2+0.5,0.3) {\small $a$};
		
		\node at (2+0.5,-0.2) {\small $b$};

	\end{tikzpicture} 
	\caption{Notations for an almost equilateral quadrilateral or simply $a^3b$-quadrilateral.}
	\label{quad}
\end{figure}
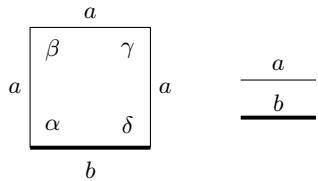

An $a^3b$-quadrilateral is given by Fig.\,\ref{quad}, with normal edge $a$, thick edge $b$ and angles $\alpha,\beta,\gamma,\delta$ as indicated. Throughout this paper, an $a^3b$-tiling is always an edge-to-edge tiling of the sphere by congruent simple quadrilaterals in Fig.\,\ref{quad}, such that all vertices have degree $\ge 3$.

Brute-force case studies of $a^3b$- or $a^4b$-tilings led to remarkable but very long papers \cite{cly,cly1}. Instead some powerful algebraic number theory approach help \cite{lw2,llhw} to classify such tilings with all angles being rational in degree (simply called \textit{rational} hereafter) very efficiently. Then the tilings with any irrational angle must have very limited vertex types and were classified easily in \cite{lw3,lsw}. 

The cost is to find all rational angle solutions to a key equation (\cite[Lemma 18]{cly}) satisfied by any  $a^3b$-quadrilateral: 
\begin{align}
	\sin\left(\aaa-\frac{\ccc}{2}\right) ~ \sin\frac{\bbb}{2}=\sin\frac{\ccc}{2} ~ \sin\left(\ddd-\frac{\bbb}{2}\right). \label{4-7}
\end{align}
Note that similar equation for $a^4b$-pentagons is much more complicated to solve in \cite{llhw}. Using $x=e^{i\theta}$ for each angle variable, \eqref{4-7} changes to polynomial form. Solving for rational angles is then equivalent to finding root of unity solutions (or cyclotomic points on some algebraic variety). According to \cite{kk}, there are two algorithms to solve such trigonometric Diophantine equations, both have super-exponential complexity-limiting practicality to either $\le12$ monomials or $\le 3$ variables. Both previous classifications \cite{lw2} and \cite{cl} applied the first algorithm,  using many old works of different authors \cite{cj,ck,mg}. Unfortunately, quite a few  typos and gaps were found in them (see Remark $17$  in \cite{lw2}, for example). 

In this paper, we will apply the other algorithm from \cite{Bradford-Davenport,BS,AS} to get an independent and also more  efficient way to find all rational $a^3b$-monotiles for the sphere, thereby supporting our full classification results more rigorously. Although the equation has $4$ variables, linear angle constraints from degree $3$ or $4$ vertices, which must appear in the tilings, can help reduce \eqref{4-7} to $2$ or $3$ variables.

	\begin{theorem*}
There are $15$ sporadic and $3$ infinite sequences of rational  $a^3b$-monotiles for the sphere, as listed in Table \ref{Tab-1.1}  and \ref{Tab-1.2}  together with all of their different tilings.
	\end{theorem*} 

\begin{table*}[htp]                        
\centering      
\begin{tabular}{c|c|c}
$f$ & $ (\aaa,\bbb,\ccc,\ddd),a,b$ & all vertices and tilings \\
\hline 
\multirow{3}{*}{$6$} 
& $(6,3,4,3)/6,1/2,1/6$ & \multirow{3}{*}{$6\aaa\bbb\ddd,2\ccc^3$ }  \\
& $(1,8,4,3)/6,0.391,1$ & \\
& $(12,4,6,2)/9,0.567,0.174$ & \\
\cline{1-3}
\multirow{5}{*}{$12$}
& $(2,10,3,6)/9,0.339,0.532$ & \multirow{5}{*}{$12\aaa\bbb\ddd,2\ccc^6$} \\				
& $(1,21,5,8)/15,0.424,0.741$ & \\
& $(4,9,5,17)/15,0.424,0.165$ & \\ 
& $(9,28,10,23)/30,0.335,0.415$ & \\ 
& $(3,16,10,41)/30,0.469,0.146$ & \\
\hline   		
\multirow{2}{*}{$20$} 
& $(5,32,6,23)/30,0.335,0.415$ & \multirow{2}{*}{$20\aaa\bbb\ddd,2\ccc^{10}$} \\
& $(1,16,6,43)/30,0.469,0.273$ & \\
\cline{1-3}		
$30$ & $(1,42,4,17)/30,0.424,0.549$ & $30\aaa\bbb\ddd,2\ccc^{15}$ \\
\hline 
\hline 
\multirow{3}{*}{$18$}
& \multirow{3}{*}{$(3,20,4,13)/18,0.339,0.452$} & $18\aaa\bbb\ddd,2\ccc^9$ \\
& & $16\aaa\bbb\ddd,2\bbb\ccc^4,2\aaa\ccc^5\ddd$ \\
& & $14\aaa\bbb\ddd,2\aaa^2\ccc\ddd^2,4\bbb\ccc^4$ \\
\hline
\hline 
$16$ & $(1,4,2,2)/4,1/4,1/2$ & $8\bbb\ddd^2,8\aaa^2\bbb\ccc,2\ccc^4$: $2$ tilings \\
\hline
$36$ & $(5,4,7,3)/9,0.174,0.258$ & $18\bbb\ccc^2,6\aaa^3\ddd,6\aaa^2\bbb^2,6\aaa\bbb\ddd^3,2\ddd^6$ \\
\hline 												
$36$ & $(15,6,10,7)/18,0.225,0.118$ & $14\aaa^2\bbb,8\aaa\ddd^3,10\bbb\ccc^3,6\bbb^2\ccc\ddd^2$ \\
\hline 
\end{tabular}
\caption{Fifteen sporadic rational $a^3b$-quadrilaterals and their tilings.}\label{Tab-1.1}        
\end{table*}

\begin{table*}[htp]               
\centering     		
\begin{tabular}{c|c}
$(\aaa,\bbb,\ccc,\ddd)$&all vertices and tilings \\
\hline 			 
\multirow{4}{*}{ $(\frac 4f,1-\frac{4}{f},\frac 4f,1)$}
& $\forall$ even $f\ge10: f\aaa\bbb\ddd,2\ccc^{\frac{f}{2}}$ \\
& $f=4k(k\ge3)$: $(f-2)\aaa\bbb\ddd,2\aaa\ccc^{\frac f4-1}\ddd,2\bbb\ccc^{\frac f4+1}$ \\
& $(f-4)\aaa\bbb\ddd,2\bbb^2\ccc^2,4\aaa\ccc^{\frac f4-1}\ddd$: $2$ tilings \\
& $f=12$: $6\aaa\bbb\ddd,2\bbb^3,6\aaa\ccc^2\ddd$ \\
\hline 			 
\multirow{4}{*}{ $(\frac 2f,\frac{4f-4}{3f},\frac{4}{f},\frac{2f-2}{3f})$}
& $\forall$ even $f\ge 6: f\aaa\bbb\ddd,2\ccc^{\frac{f}{2}}$ \\
& $f=6k+4(k\ge1)$: $(f-2)\aaa\bbb\ddd,2\bbb\ccc^{\frac{f+2}{6}},2\aaa\ccc^{\frac{f-1}{3}}\ddd$ \\
& $(f-4)\aaa\bbb\ddd,2\aaa^2\ccc^{\frac{f-4}{6}}\ddd^2,4\bbb\ccc^{\frac{f+2}{6}}$: $\lfloor \frac{k+2}{2} \rfloor$ tilings \\
& $(f-6)\aaa\bbb\ddd,2\aaa\ddd^3,2\aaa^2\bbb\ccc^{\frac{f-4}{6}},4\bbb\ccc^{\frac{f+2}{6}}$: $3$ tilings \\
\hline 
\multirow{4}{*}{ $(\frac2f,\frac{2f-4}{3f},\frac4f,\frac{4f-2}{3f})$}
& $\forall$ even $f\ge10: f\aaa\bbb\ddd,2\ccc^{\frac{f}{2}}$ \\
& $f=6k+2(k\ge2)$: $(f-2) \aaa\bbb\ddd,2\aaa\ccc^{\frac{f-2}{6}}\ddd,2\bbb\ccc^{\frac{f+1}{3}}$ \\
& $(f-4)\aaa\bbb\ddd,4\aaa\ccc^{\frac{f-2}{6}}\ddd,2\bbb^2\ccc^{\frac{f+4}{6}}$: $\lfloor \frac{k+3}{2} \rfloor$ tilings \\
& $(f-6)\aaa\bbb\ddd,2\bbb^3\ccc,6\aaa\ccc^{\frac{f-2}{6}}\ddd$ \\
\hline
\end{tabular}
\caption{Three infinite sequences of  rational $a^3b$-quadrilaterals and their tilings.}\label{Tab-1.2}        
\end{table*}
 
 In Table \ref{Tab-1.1} and \ref{Tab-1.2}  the angles and edge lengths are expressed in units of $\pi$,  and the last column counts all vertices and also all tilings when they are not uniquely determined by the vertices. A rational fraction, such as $\alpha=\frac{1}{6}$, means the precise value $\frac{\pi}{6}$. A decimal expression, such as $a\approx0.391$, means an approximate value  $0.391\pi \le a < 0.392\pi$. 
 
 In conclusion, this paper verifies that the list of rational  $a^3b$-monotiles in our previous work \cite{lw2} was complete while it was incomplete in \cite{cl}. 
 
\subsection*{Outline of the paper}

This paper is organized as follows. Section \ref{basic_facts} includes basic facts from \cite{lw1,lw3} and some technical results specific to $a^3b$. Section \ref{region}  introduces the method of solving algebraic equations in roots of unity with some examples. Section \ref{cases} applies this method to all possible degree 3 vertex types, completing the proof.

	\section{Basic Facts}
	\label{basic_facts}
	We will always express angles in $\pi$ radians for simplicity. So the sum of all angles at a vertex is $2$. We present some basic facts and techniques in this section. 

    Let $v,e,f$ be the numbers of vertices, edges, and tiles in a quadrilateral tiling. Let $v_k$ be the number of vertices of degree $k$. Euler's formula $v-e+f=2$ implies (see \cite{lw1})
\begin{align}
	f&=6+ \sum_{k=4}^{\infty}(k-3)v_k
	=6+v_4+2v_5+3v_6+\cdots, \label{vcountf} \\
	v_3 &=8+\sum_{k=5}^{\infty}(k-4)v_k=8+v_5+2v_6+3v_7+\cdots. \label{vcountv}
\end{align}
    So $f\ge 6$ and $v_3 \ge 8$. 

\begin{lemma} [{\cite[Lemma 2]{lw1}}] \label{anglesum} 
	If all tiles in a tiling of the sphere by $f$ quadrilaterals have the same four angles $\aaa,\bbb,\ccc,\ddd$, then 
	\[
	\aaa+\bbb+\ccc+\ddd = 2+\frac{4}{f} , 
	\]
	ranging in $(2,\frac83]$. In particular no vertex contains all four angles.
\end{lemma}

\begin{lemma} [{\cite[Lemma 3]{wy2}}] \label{geometry1}
	If the  quadrilateral in Fig.\,\ref{quad} is simple, then $\bbb<\ccc$ is equivalent to $\aaa>\ddd$.
\end{lemma}

\begin{lemma}[{\cite[Lemma 3]{lw2}}]\label{geometry3}
	If the quadrilateral in Fig.\,\ref{quad} is simple, then $\bbb=\ddd$ if and only if $\aaa=1$. Furthermore, if it is convex with all angles $<1$, then $\bbb>\ddd$ is equivalent to $\aaa<\ccc$, and  $\bbb<\ddd$ is equivalent to $\aaa>\ccc$.    
\end{lemma}

\begin{lemma}[{\cite[Lemma 4]{lw2}}]\label{geometry4}
	If the quadrilateral in Fig.\,\ref{quad} is simple, and $\ddd\le 1$, then $2\aaa+\bbb>1$ and $\bbb+2\ccc>1$ .
\end{lemma}

\begin{lemma}[{\cite[Lemma 8]{lw2}}]\label{proposition-1}
	There is no tiling of the sphere by congruent quadrilaterals with two angles $\ge1$. 
\end{lemma}

\begin{lemma}[{\cite[Lemma 2.9]{lw3}}]\label{relation}
	Any convex spherical $a^3b$-quadrangle with all angles $<1$ satisfies 
	\[ 	\aaa+\ddd<1+\bbb, \quad	\aaa+\ddd<1+\ccc, \quad	\aaa+\bbb<1+\ddd, \quad	\ccc+\ddd<1+\aaa.\]
\end{lemma}

\begin{lemma} [{Parity Lemma, \cite[Lemma 10]{wy2}}] \label{beven}
	In an $a^3b$-tiling, the total number of $ab$-angles $\aaa$ and $\ddd$ at any vertex is even. 
\end{lemma}

\begin{lemma}[{Balance Lemma, \cite[Lemma 11]{wy2}}] \label{balance}
	 In an $a^3b$-tiling by $f$ congruent tiles, each angle of the tile appears $f$ times in total. If either $\aaa^2\cdots$ or $\delta^2\cdots$ is not a vertex, then any vertex either has no $\aaa,\delta$, or is of the form $\aaa\delta\cdots$ with no more $\aaa,\delta$ in the remainder.
\end{lemma}

\begin{lemma} [{\cite[Proposition 10]{lw2}}] \label{symmetric}	
	There is no tiling of the sphere by congruent symmetric $a^3b$-quadrilaterals ($\aaa=\ddd$ and $\bbb=\ccc$).
\end{lemma}


In an $a^3b$-tiling, there are only ten possible degree $3$ vertices:
four $a^3$-vertices $\bbb^3$, $\bbb^2\ccc$, $\bbb\ccc^2$, $\ccc^3$ and six $a^2b$-vertices $\aaa^2\bbb$, $\aaa^2\ccc$, $\bbb\ddd^2$, $\ccc\ddd^2$, $\aaa\bbb\ddd$, $\aaa\ccc\ddd$, 

\begin{lemma}[{\cite[Lemma 2.11]{lw3}}]\label{geometry2}
	In an $a^3b$-tiling, we cannot have three different $a^2b$-vertices. Moreover, the following are the only possible combinations of two $a^2b$-vertices:
	\begin{enumerate}
		\item $\aaa\bbb\delta,\ccc\delta^2$.
		\item $\aaa\ccc\delta,\aaa^2\bbb$.
		\item $\aaa^2\bbb,\ccc\delta^2$.
	\end{enumerate}
\end{lemma}

When $\aaa\bbb\ddd$ (or $\aaa\ccc\ddd$) appears as a vertex, we get $\ccc=\frac{4}{f}$ and such a quadrilateral always admits the simplest $2$-layer earth map tiling. When $\aaa\bbb\ddd$ (or $\aaa\ccc\ddd$) never appears as a vertex, we have the following two Lemmas.

\begin{lemma} [{\cite[Table 5]{lw3}}]\label{tildev20}
	In an $a^3b$-tiling, if there are two (or more) different degree $3$ vertex types without $\aaa\bbb\ddd$, $\aaa\ccc\ddd$, then they must be one of the following: $\{\aaa^2\bbb, \bbb^3\}$, $\{\aaa^2\bbb, \bbb^2\ccc\}$, $\{\aaa^2\bbb, \ccc^3\}$, $\{\aaa^2\bbb, \ccc\ddd^2\}$, $\{\bbb^2\ccc, \bbb\ddd^2\}$, $\{\bbb\ccc^2, \bbb\ddd^2\}$ or $\{\bbb\ddd^2, \ccc^3\}$, up to the symmetry of interchanging $\aaa\leftrightarrow\ddd$ and $\bbb\leftrightarrow\ccc$.
\end{lemma}

If there is a unique degree $3$ vertex type without $\aaa\bbb\ddd$, $\aaa\ccc\ddd$, then it must be $\aaa^2\bbb,\bbb\ddd^2,\bbb^3$ or $\bbb\ccc^2$, up to the symmetry of interchanging $\aaa\leftrightarrow\ddd$ and $\bbb\leftrightarrow\ccc$.

\begin{lemma} [{\cite[Lemma 5.1 and Lemma 5.2]{lw3}}]\label{combined_lemma}
	In an $a^3b$-tiling, we have the following:
	\begin{enumerate}
		\item If all degree $3$ vertices are of type $\bbb\ccc^2$ or $\bbb^3$, then there exists some degree $4$ vertex with two $b$-edges: one of ${\aaa^4,\aaa^3\ddd,\aaa^2\ddd^2,\aaa\ddd^3,\ddd^4}$ must appear.
		\item If all degree $3$ vertices are of type $\aaa^2\bbb$, then there exists a degree $4$ vertex without $\aaa$: one of ${\bbb^4,\bbb^3\ccc,\bbb^2\ccc^2,\bbb\ccc^3,\bbb^2\ddd^2,\bbb\ccc\ddd^2,\ccc^4,\ccc^2\ddd^2,\ddd^4}$ must appear.
		\item If all degree $3$ vertices are of type $\bbb\ddd^2$, then there exists a degree $4$ vertex without $\ddd$: one of ${\aaa^4,\aaa^2\bbb^2,\aaa^2\ccc^2,\aaa^2\bbb\ccc,\bbb^4,\bbb^3\ccc,\bbb^2\ccc^2,\bbb\ccc^3,\ccc^4}$ must appear.
	\end{enumerate}
\end{lemma}

\section{An explicit example applying the number theory techniques}
\label{region}

To find all rational $a^3b$-monotiles, we will find all reasonable rational angle solutions to \eqref{4-7} under $36$ different linear angle constraints in the next section. The computations can be done by SageMath, Maple or other symbolic computation software.

Take Case $\{\bbb^3,\aaa^4\}$ in Table \ref{3ab2} as a typical example, and we will show explicitly how to resolve it by the number theory techniques developed by Bradford-Davenport \cite{Bradford-Davenport} and Beukers--Smyth \cite{BS}. 

Firstly by solving the linear system of angle sums, we get
\[
\aaa=\frac{\pi}{2}, \quad \bbb=\frac{2\pi}{3},\quad \ccc=\frac{5\pi}{6}-\ddd+\frac{4\pi}{f}. \]
Then the trigonometric equation \eqref{4-7}  reduces to 
\begin{align*}
	Q(\delta, f) &:= 2 \sin(\frac{\pi}{4}+\frac{\delta}{2}-\frac{2 \pi}{f})-2 \cos(\frac{5 \pi}{12}+\frac{\delta}{2}-\frac{2 \pi}{f})+2 \cos(\frac{\pi}{12}+\frac{\delta}{2}+\frac{2 \pi}{f})+2 \cos(\frac{\pi}{4}+\frac{3 \delta}{2}-\frac{2 \pi}{f})=0	
\end{align*}
By setting $x={\mathrm e}^{\mathrm{i} \ddd},y={\mathrm e}^{\frac{4\mathrm{i} \pi}f}$,
the above equation is transformed to a polynomial one: 
\begin{equation*}
	P(x,y) := \zeta_{12}^{4} x^{3}+\zeta_{12}^{3} x^{2} y-(\zeta_{12}^{5}-\zeta_{12}) x^{2}+(\zeta_{12}^{4}-1) x y+\zeta_{12}^{2} x+\zeta_{12} y = 0.	
\end{equation*}
Replace $\zeta_{12}$ by the other primitive $12$-th roots of unity $\zeta_{12}^k$ for $k=5,7,11$ in $P(x,y)$, then the product of these four polynomials will have rational coefficients by Galois theory: 
$\tilde{P}(x,y)=x^{12}-x^{10} y^{2}+x^{8} y^{4}+12 x^{10} y+5 x^{10}+43 x^{8} y^{2}+5 x^{6} y^{4}+12 x^{8} y+48 x^{6} y^{3}+24 x^{8}+24 x^{6} y^{2}+24 x^{4} y^{4}+48 x^{6} y+12 x^{4} y^{3}+5 x^{6}+43 x^{4} y^{2}+5 x^{2} y^{4}+12 x^{2} y^{3}+x^{4}-x^{2} y^{2}+y^{4}.$

Note that setting $x={\mathrm e}^{2\mathrm{i} \ddd},y={\mathrm e}^{\frac{4\mathrm{i} \pi}f}$ gives a simpler polynomial $L(x,y) := x^{4} y^{4}-x^{5} y^{2}+5 x^{3} y^{4}+x^{6}+12 x^{5} y+43 x^{4} y^{2}+48 x^{3} y^{3}+24 x^{2} y^{4}+5 x^{5}+12 x^{4} y+24 x^{3} y^{2}+12 x^{2} y^{3}+5 x \,y^{4}+24 x^{4}+48 x^{3} y+43 x^{2} y^{2}+12 x \,y^{3}+y^{4}+5 x^{3}-x \,y^{2}+x^{2}$, which is full (i.e. its monomials' exponent vectors in $\Z^2$ generate the full integer lattice $\Z^2$). 

To find roots of unity solutions $(x,y)$, or to find all cyclotomic points on this algebraic curve, we follow \cite{BS} to compute the resultants of $L(x,y)$ with $L(-x,y)$, $L(x,-y)$, $L(-x,-y)$, $L(x^2,y^2)$, $L(-x^2,y^2)$, $L(x^2,-y^2)$, and $L(-x^2,-y^2)$ respectively. For example, 
$$
\begin{array}{lcc}
	f(x) :=& {\mathrm{Resultant} (\mathit{L(x,y)} ,\mathit{L(-x,y)} ,x)\, }  =  256 y^{8} 	(25 y^{8}+998 y^{6}+1923 y^{4}+998 y^{2}+25)^{2}  \\
	 & (11 y^{8}+6 y^{7}-1217 y^{6}-4182 y^{5}-5112 y^{4}-4182 y^{3}-1217 y^{2}+6 y+11)^{2}.  
\end{array}
$$

This is a single-variable polynomial, whose root of unity solutions (or cyclotomic part) can be found by a recursive algorithm \cite{Bradford-Davenport} computing $\gcd(f(x), f(-x))$, $\gcd(f(x), f(x^2))$, and $\gcd(f(x), f(-x^2))$.  Both \cite{Bradford-Davenport} and  \cite{BS} are based on a key but simple  property: For $\omega$ a root of unity of order $k \equiv r \pmod{4}$ ($r=0,1,2,3$), $\omega$ is conjugate to $-\omega$, $\omega^2$, $-\omega^2$, $\omega^2$ respectively. It turns out that all such $\gcd$ for $f(x)$ are $1$, thus $f(x)$ has no root of unity solution.

Similarly, the following resultant produces two root of unity solutions $y=\zeta_4^i$ for $i = 1, 3$: 
$$
\begin{array}{lcc}
	{\mathrm{Resultant} (\mathit{L(x,y)} ,\mathit{L(x,-y)} ,x)\, }  &  \\
	= 247669456896 y^{14} (4 y^{8}+67 y^{6}+114 y^{4}+67 y^{2}+4) (y^{8}+72 y^{6}+110 y^{4}+72 y^{2}+1) (y^{2}+1)^{2}.  &
\end{array}
$$

After all seven resultants computations, we get totally six root of unity solutions $y=\zeta_{12}^i$ for $i = 2,3,4,8,9,10$. By $y = e^{\frac{4i\pi}{f}}$, we get $f=12$, $8$, $6$, $3$, $\frac{8}{3}$, and $\frac{12}{5}$. Since $f$ must be an even integer $\ge6$, we only need to consider $i = 2,3,4$. 

For $f=12$, substitute $y = e^{\frac{i\pi}{3}}$ in $P(x,y)=0$, and we get similarly $x=\zeta_{12}^i$ for $i = 2,3,4,8,9,10$. This implies $(\aaa,\bbb,\ccc,\ddd)=(3, 4, 2, 5)/6$, $(3, 4, 3, 4)/6$, $(3, 4, 4, 3)/6$, $(3, 4, -4, 11)/6$, $(3, 4, -3, 10)/6$ or $(3, 4, -2, 9)/6$. The last three solutions are dismissed since they have negative angles. The first solution is dismissed since it does not satisfy \eqref{4-7} and is a fake solution. The second and third solutions are dismissed since they violate Lemma \ref{geometry3} and \ref{symmetric} respectively.

Similarly, for $f=6$ we get $(\aaa,\bbb,\ccc,\ddd)=(3,4,8,1)/6$, and for $f=8$ we get $(\aaa,\bbb,\ccc,\ddd)=(3,4,6,2)/6$. Only these two good solutions satisfying  \eqref{4-7} and Lemma \ref{anglesum},\ref{geometry1},\ref{geometry3},\ref{geometry4},\ref{proposition-1},\ref{relation}, are collected into Table \ref{3ab2}. 

It turns out that $(\aaa,\bbb,\ccc,\ddd)=(3,4,8,1)/6$ is a monotile admitting the standard 2-layer earth map tiling and it also appears in Table \ref{abd} after the symmetry $\aaa\leftrightarrow\ddd$, $\bbb\leftrightarrow\ccc$. 

However, $(\aaa,\bbb,\ccc,\ddd)=(3,4,6,2)/6$ does not admit any tiling of the sphere. Otherwise the parity lemma implies that $\ccc\cdots=\aaa^2\ccc$. Consequently, we have $\#\aaa=2\#\ccc=2f$, a contradiction to the balance lemma.  Similarly, the following two cases in Table \ref{3ab2} produce no monotiles either: 
\begin{enumerate}
	\item Case $\{\aaa^2\bbb,\bbb^4,(9,6,12,5)/12,6\}$.
	\item Case $\{\bbb\ddd^2,\aaa^2\ccc^2,(2,6,4,3)/6,8\}$.
\end{enumerate}

The above discussion resolves Case $\{\bbb^3,\aaa^4\}$ completely. For more details or more examples, we refer to our work \cite{llhw} on rational $a^4b$-tilings.

\section{All rational $a^3b$-monotiles}\label{cases}

 By Lemma \ref{symmetric}, we only need to consider non-symmetric $a^3b$-quadrilaterals ($\aaa\neq\ddd$, $\bbb\neq\ccc$) with even $f \geq 6$  and all rational angles in $(0,2\pi)$ satisfying \eqref{4-7} and Lemma \ref{anglesum},\ref{geometry1},\ref{geometry3},\ref{geometry4},\ref{proposition-1},\ref{relation} (simply called ``good'' hereafter). By Lemma \ref{geometry2},\ref{tildev20},\ref{combined_lemma}, there are totally $36$ vertex combinations to help simplify \eqref{4-7}, which induce all rational $a^3b$-monotiles summarized in Tables  \ref{abd},\ref{3ab1} and \ref{3ab2} of the following subsections respectively.

\subsection{Tilings with $\aaa\bbb\ddd$ or $\aaa\ccc\ddd$ as a vertex}

By the symmetry, it is enough to consider $\aaa\bbb\ddd$.
By Lemma \ref{anglesum}, we get
\[
\aaa=2\pi-\bbb-\ddd, \quad \ccc=\tfrac{4\pi}{f}. 
\]	
The trigonometric equation \eqref{4-7}  reduces to
\begin{align*}
	Q&(\bbb, \ddd, f):= 2 \cos \! \left(\tfrac{3 \beta}{2}+\delta +\tfrac{2 \pi}{f}\right)-2 \cos \! \left(\tfrac{\beta}{2}+\delta +\tfrac{2 \pi}{f}\right)+2 \cos \! \left(\tfrac{\beta}{2}-\tfrac{2 \pi}{f}-\delta \right)-2 \cos \! \left(\tfrac{\beta}{2}-\delta +\tfrac{2 \pi}{f}\right)=0
\end{align*}
By setting $x={\mathrm e}^{\frac{\mathrm{i\bbb}}2},y={\mathrm e}^{\mathrm{2i} \delta},z={\mathrm e}^{\frac{4\mathrm{i} \pi}f}$, $Q(\bbb,\ddd,f)$  is transformed to a polynomial
\begin{align*}
	P(x,y,z)&:=x^{3} y z -x^{2}y z-x^{2} z +x y z +x^{2}-x y -x +1=0.
\end{align*}
This is  the only case inducing a $3$-variable polynomial equation.
Following the algorithm in \cite{AS}, we compute $15$ resultants to get all good rational angle solutions below:

\subsubsection*{Case 1. resultant$(P(x,y,z),P(-x,y,z),z)$}

$\left(z -1\right)^{2} \left(y +1\right)^{2} \left(y z -1\right)^{2}=0$.

\begin{table}[H]        
	\centering 
	\begin{tabular}{|c|c|c|c|}
		\hline  
		$f$& $(\aaa,\bbb,\ccc,\ddd)$&$f$& $(\aaa,\bbb,\ccc,\ddd)$ \\
		\hline\hline
		$6$&$(6,3,4,3)/6$&$6$&$(1,8,4,3)/6$\\
		\hline  
	\end{tabular}
\end{table}

\subsubsection*{Case 2. resultant$(P(x,y,z),P(x,-y,z),z)$}

$\left(z^{2}+1\right) \left(z -1\right)^{2}=0$.

\begin{table}[H]        
	\centering 
	\begin{tabular}{|c|c|}
		\hline  
		$f$& $(\aaa,\bbb,\ccc,\ddd)$\\
		\hline\hline
		$8$&$(3,14,6,7)/12$\\
		\hline  
	\end{tabular}
\end{table}
\subsubsection*{Case 3. resultant$(P(x,y,z),P(x,y,-z),z)$}

$\left(y -1\right)^{2} \left(y +1\right)^{2}=0$

\begin{table}[H]        
	\centering 
	\begin{tabular}{|c|c|c|c|c|c|}
		\hline  
		$f$& $(\aaa,\bbb,\ccc,\ddd)$&$f$& $(\aaa,\bbb,\ccc,\ddd)$&$f$& $(\aaa,\bbb,\ccc,\ddd)$ \\
		\hline\hline
		$6$&$(6,3,4,3)/6$&$6$&$(1,8,4,3)/6$&&$(4,f-4,4,f)/f$\\
		\hline  
	\end{tabular}
\end{table}
\subsubsection*{Case 4.  resultant$(P(\tilde{x},\tilde{y},z),P(-\tilde{x},-\tilde{y},z),z), (x=\tilde{x}\tilde{y}, y=\tfrac{\tilde{x}}{\tilde{y}})$}

$\left(y -1\right) \left(y +1\right)=0$

\begin{table}[H]        
	\centering 
	\begin{tabular}{|c|c|c|c|c|c|}
		\hline  
		$f$& $(\aaa,\bbb,\ccc,\ddd)$&$f$& $(\aaa,\bbb,\ccc,\ddd)$&$f$& $(\aaa,\bbb,\ccc,\ddd)$ \\
		\hline\hline
		$6$&$(6,3,4,3)/6$&$6$&$(1,8,4,3)/6$&&$(4,f-4,4,f)/f$\\
		\hline  
	\end{tabular}
\end{table}
\subsubsection*{Case 5.  resultant$(P(\tilde{x},y,z),P(-\tilde{x},y,-z),z), (x=\tilde{x}^2)$}

$\left(y +1\right)^{3} \left(z^{2}+1\right) \left(y -1\right)=0$

\begin{table}[H]        
	\centering 
	\begin{tabular}{|c|c|c|c|c|c|c|c|}
		\hline  
		$f$& $(\aaa,\bbb,\ccc,\ddd)$&$f$& $(\aaa,\bbb,\ccc,\ddd)$&$f$& $(\aaa,\bbb,\ccc,\ddd)$&$f$& $(\aaa,\bbb,\ccc,\ddd)$ \\
		\hline\hline
		$6$&$(6,3,4,3)/6$&$6$&$(1,8,4,3)/6$&$8$&$(3,14,6,7)/12$&&$(4,f-4,4,f)/f$\\
		\hline  
	\end{tabular}
\end{table}

\subsubsection*{Case 6.  resultant$(P(x,y,z),P(x,-y,-z),z)$}

$\left(y -1\right) \left(y +1\right) \left(y^{2}+y z +z^{2}\right)=0$

\begin{table}[H]        
	\centering 
	\begin{tabular}{|c|c|c|c|c|c|c|c|}
		\hline  
		$f$& $(\aaa,\bbb,\ccc,\ddd)$&$f$& $(\aaa,\bbb,\ccc,\ddd)$&$f$& $(\aaa,\bbb,\ccc,\ddd)$&$f$& $(\aaa,\bbb,\ccc,\ddd)$ \\
		\hline\hline
		$6$&$(6,3,4,3)/6$&$6$&$(1,8,4,3)/6$&$8$&$(3,14,6,7)/12$&&$(4,f-4,4,f)/f$\\
		\hline  
	\end{tabular}
\end{table}

\subsubsection*{Case 7.  resultant$(P(\tilde{x},\tilde{y},z),P(-\tilde{x},-\tilde{y},-z),z), (x=\tilde{x}\tilde{y}, y=\tfrac{\tilde{x}}{\tilde{y}})$}

$\left(z^{2}+1\right) \left(y^{3} z -1\right) \left(y z -1\right)=0$

\begin{table}[H]        
	\centering 
	\begin{tabular}{|c|c|c|c|c|c|c|c|}
		\hline  
		$f$& $(\aaa,\bbb,\ccc,\ddd)$&$f$& $(\aaa,\bbb,\ccc,\ddd)$&$f$& $(\aaa,\bbb,\ccc,\ddd)$&$f$& $(\aaa,\bbb,\ccc,\ddd)$ \\
		\hline\hline
		$6$&$(12,4,6,2)/9$&$8$&$(3,14,6,7)/12$&&$(6,4f-4,12,2f-2)/f$&&$(6,2f-4,12,4f-2)/f$\\
		\hline  
	\end{tabular}
\end{table}

\subsubsection*{Case 8.  resultant$(P(x,y,z),P(x^2,y^2,z^2),z)$}

$\left(z -1\right)^{4} \left(y z -1\right)^{3} \left(y -1\right)^{2} \left(y^{3} z -1\right)=0$

\begin{table}[H]        
	\centering 
	\begin{tabular}{|c|c|c|c|c|c|c|c|}
		\hline  
		$f$& $(\aaa,\bbb,\ccc,\ddd)$&$f$& $(\aaa,\bbb,\ccc,\ddd)$&$f$& $(\aaa,\bbb,\ccc,\ddd)$&$f$& $(\aaa,\bbb,\ccc,\ddd)$ \\
		\hline\hline
		$6$&$(12,4,6,2)/9$&&$(4,f-4,4,f)/f$&&$(6,4f-4,12,2f-2)/f$&&$(6,2f-4,12,4f-2)/f$\\
		\hline  
	\end{tabular}
\end{table}
\subsubsection*{Case 9.  resultant$(P(x,y,z),P(-x^2,y^2,z^2),z)$}

$\left(y +1\right)^{3} \left(z -1\right)^{2} \left(y -1\right) \left(y^{2}+y z +z^{2}\right) \left(y z -1\right)^{2}=0$

\begin{table}[H]        
	\centering 
	\begin{tabular}{|c|c|c|c|c|c|c|c|}
		\hline  
		$f$& $(\aaa,\bbb,\ccc,\ddd)$&$f$& $(\aaa,\bbb,\ccc,\ddd)$&$f$& $(\aaa,\bbb,\ccc,\ddd)$&$f$& $(\aaa,\bbb,\ccc,\ddd)$ \\
		\hline\hline
		$6$&$(6,3,4,3)/6$&$6$&$(1,8,4,3)/6$&$8$&$(3,14,6,7)/12$&&$(4,f-4,4,f)/f$\\
		\hline  
	\end{tabular}
\end{table}
\subsubsection*{Case 10.  resultant$(P(x,y,z),P(x^2,-y^2,z^2),z)$}

$(z -1)^{2}(y -1) (y +1) (y^{6} z^{4}-2 y^{6} z^{3}-2 y^{5} z^{4}+y^{6} z^{2}-y^{4} z^{4}+4 y^{5} z^{2}+7 y^{4} z^{3}+2 y^{3} z^{4}-2 y^{5} z -y^{4} z^{2}+y^{3} z^{3}+y^{2} z^{4}-5 y^{4} z -10 y^{3} z^{2}-5 y^{2} z^{3}+y^{4}+y^{3} z -y^{2} z^{2}-2 yz^{3}+2 y^{3}+7 y^{2} z +4 yz^{2}-y^{2}+z^{2}-2 y -2 z +1)=0$

\begin{table}[H]        
	\centering 
	\begin{tabular}{|c|c|c|c|c|c|c|c|}
		\hline  
		$f$& $(\aaa,\bbb,\ccc,\ddd)$&$f$& $(\aaa,\bbb,\ccc,\ddd)$&$f$& $(\aaa,\bbb,\ccc,\ddd)$&$f$& $(\aaa,\bbb,\ccc,\ddd)$ \\
		\hline\hline
		$6$&$(6,3,4,3)/6$&$6$&$(1,8,4,3)/6$&$8$&$(3,14,6,7)/12$&$16$&$(1,3,1,4)/4$\\
		\hline 
		$16$&$(1,10,2,5)/8$&$16$&$(3,14,6,31)/24$&&$(4,f-4,4,f)/f$&&\\ 
		\hline
	\end{tabular}
\end{table}

\subsubsection*{Case 11.  resultant$(P(x,y,z),P(x^2,y^2,-z^2),z)$}

$ (y -1)^3 (y +1)(y^{4} z^{6}-3 y^{4} z^{5}+4 y^{4} z^{4}-4 y^{3} z^{5}-2 y^{4} z^{3}+10 y^{3} z^{4}-y^{2} z^{5}+y^{4} z^{2}-9 y^{3} z^{3}+8 y^{2} z^{4}-yz^{5}+4 y^{3} z^{2}-12 y^{2} z^{3}+4 yz^{4}-y^{3} z +8 y^{2} z^{2}-9 yz^{3}+z^{4}-y^{2} z +10 yz^{2}-2 z^{3}-4 y z +4 z^{2}-3 z +1)=0$

\begin{table}[H]        
	\centering 
	\begin{tabular}{|c|c|c|c|c|c|c|c|}
		\hline  
		$f$& $(\aaa,\bbb,\ccc,\ddd)$&$f$& $(\aaa,\bbb,\ccc,\ddd)$&$f$& $(\aaa,\bbb,\ccc,\ddd)$&$f$& $(\aaa,\bbb,\ccc,\ddd)$ \\
		\hline\hline
		$6$&$(6,3,4,3)/6$&$6$&$(1,8,4,3)/6$&$8$&$(3,14,6,7)/12$&$12$&$(2,10,3,6)/9$\\
		\hline 
		$12$&$(1,21,5,8)/15$&$12$&$(4,9,5,17)/15$&$12$&$(9,28,10,23)/30$&$12$&$(3,16,10,41)/30$\\ 
		\hline
		$12$&$(1,2,1,3)/3$&$12$&$(3,22,6,11)/18$&$12$&$(3,10,6,23)/18$&$24$&$(1,5,1,6)/6$\\
		\hline
		$24$&$(3,46,6,23)/36$&$24$&$(3,22,6,47)/36$&$36$&$(1,8,1,9)/9$&$36$&$(3,70,6,35)/54$\\
		\hline
		$36$&$(3,37,6,71)/54$&&$(4,f-4,4,f)/f$&&&&\\
		\hline
	\end{tabular}
\end{table}

\subsubsection*{Case 12.   resultant$(P(x,y,z),P(-x^2,-y^2,z^2),z)$}

$ y^{8} z^{6}+y^{7} z^{7}+y^{6} z^{8}-2 y^{8} z^{5}-y^{6} z^{7}+y^{8} z^{4}-4 y^{7} z^{5}-y^{6} z^{6}+6 y^{7} z^{4}+3 y^{6} z^{5}+y^{5} z^{6}+y^{4} z^{7}-3 y^{7} z^{3}+2 y^{6} z^{4}+3 y^{5} z^{5}-3 y^{4} z^{6}-y^{3} z^{7}-6 y^{6} z^{3}-6 y^{5} z^{4}+5 y^{4} z^{5}+y^{3} z^{6}+3 y^{6} z^{2}+4 y^{5} z^{3}-4 y^{4} z^{4}+4 y^{3} z^{5}+3 y^{2} z^{6}+y^{5} z^{2}+5 y^{4} z^{3}-6 y^{3} z^{4}-6 y^{2} z^{5}-y^{5} z -3 y^{4} z^{2}+3 y^{3} z^{3}+2 y^{2} z^{4}-3 yz^{5}+y^{4} z +y^{3} z^{2}+3 y^{2} z^{3}+6 yz^{4}-y^{2} z^{2}-4 yz^{3}+z^{4}-y^{2} z -2 z^{3}+y^{2}+y z +z^{2}=0$

\begin{table}[H]        
	\centering 
	\begin{tabular}{|c|c|c|c|c|c|c|c|}
		\hline  
		$f$& $(\aaa,\bbb,\ccc,\ddd)$&$f$& $(\aaa,\bbb,\ccc,\ddd)$&$f$& $(\aaa,\bbb,\ccc,\ddd)$&$f$& $(\aaa,\bbb,\ccc,\ddd)$ \\
		\hline\hline
		$8$&$(3,14,6,7)/12$&$16$&$(1,3,1,4)/4$&$16$&$(1,10,2,5)/8$&$16$&$(3,14,6,31)/24$\\
		\hline 
		$30$&$(1,42,4,17)/30$&$30$&$(2,13,2,15)/15$&$30$&$(3,58,6,29)/45$&$30$&$(3,28,6,59)/45$\\ 
		\hline
		$42$&$(2,19,2,21)/21$&$42$&$(3,82,6,41)/63$&$42$&$(3,40,6,83)/63$&&\\
		\hline
	\end{tabular}
\end{table}

\subsubsection*{Case 13.  resultant$(P(x,y,z),P(-x^2,y^2,-z^2),z)$}

$(z^{2}+1)(y +1)^{2} (y^{6} z^{4}+y^{5} z^{5}+y^{4} z^{6}-y^{6} z^{3}-3 y^{5} z^{4}-2 y^{4} z^{5}+y^{6} z^{2}+3 y^{5} z^{3}+3 y^{4} z^{4}-y^{3} z^{5}-y^{5} z^{2}+y^{4} z^{3}+4 y^{3} z^{4}+2 y^{2} z^{5}-2 y^{4} z^{2}-6 y^{3} z^{3}-2 y^{2} z^{4}+2 y^{4} z +4 y^{3} z^{2}+y^{2} z^{3}-yz^{4}-y^{3} z +3 y^{2} z^{2}+3 yz^{3}+z^{4}-2 y^{2} z -3 yz^{2}-z^{3}+y^{2}+y z +z^{2})=0$

\begin{table}[H]        
	\centering 
	\begin{tabular}{|c|c|c|c|c|c|c|c|}
		\hline  
		$f$& $(\aaa,\bbb,\ccc,\ddd)$&$f$& $(\aaa,\bbb,\ccc,\ddd)$&$f$& $(\aaa,\bbb,\ccc,\ddd)$&$f$& $(\aaa,\bbb,\ccc,\ddd)$ \\
		\hline\hline
		$6$&$(6,3,4,3)/6$&$6$&$(1,8,4,3)/6$&$8$&$(3,14,6,7)/12$&$12$&$(2,10,3,6)/9$\\
		\hline 
		$12$&$(1,21,5,8)/15$&$12$&$(4,9,5,17)/15$&$12$&$(9,28,10,23)/30$&$12$&$(3,16,10,41)/30$\\ 
		\hline
		$12$&$(1,2,1,3)/3$&$12$&$(3,22,6,11)/18$&$12$&$(3,10,6,23)/18$&$16$&$(1,3,1,4)/4$\\
		\hline
		$16$&$(1,10,2,5)/8$&$16$&$(3,14,6,31)/24$&$60$&$(1,14,1,15)/15$&$60$&$(3,118,6,59)/90$\\
		\hline
		$60$&$(3,58,6,119)/90$&&&&&&\\
		\hline
	\end{tabular}
\end{table}
\subsubsection*{Case 14.  resultant$(P(x,y,z),P(x^2,-y^2,-z^2),z)$}

$y^{8} z^{8}-4 y^{8} z^{7}-2 y^{7} z^{8}+8 y^{8} z^{6}+5 y^{7} z^{7}-y^{6} z^{8}-9 y^{8} z^{5}-5 y^{7} z^{6}+9 y^{6} z^{7}+2 y^{5} z^{8}+7 y^{8} z^{4}-24 y^{6} z^{6}-8 y^{5} z^{7}+y^{4} z^{8}-3 y^{8} z^{3}+4 y^{7} z^{4}+35 y^{6} z^{5}+12 y^{5} z^{6}-7 y^{4} z^{7}+y^{8} z^{2}-3 y^{7} z^{3}-30 y^{6} z^{4}-5 y^{5} z^{5}+23 y^{4} z^{6}+3 y^{3} z^{7}+y^{7} z^{2}+18 y^{6} z^{3}-4 y^{5} z^{4}-41 y^{4} z^{5}-6 y^{3} z^{6}+2 y^{2} z^{7}-8 y^{6} z^{2}+8 y^{5} z^{3}+48 y^{4} z^{4}+8 y^{3} z^{5}-8 y^{2} z^{6}+2 y^{6} z -6 y^{5} z^{2}-41 y^{4} z^{3}-4 y^{3} z^{4}+18 y^{2} z^{5}+y z^{6}+3 y^{5} z +23 y^{4} z^{2}-5 y^{3} z^{3}-30 y^{2} z^{4}-3 y z^{5}+z^{6}-7 y^{4} z +12 y^{3} z^{2}+35 y^{2} z^{3}+4 y z^{4}-3 z^{5}+y^{4}-8 y^{3} z -24 y^{2} z^{2}+7 z^{4}+2 y^{3}+9 y^{2} z -5 y z^{2}-9 z^{3}-y^{2}+5 y z +8 z^{2}-2 y -4 z +1=0$

\begin{table}[H]        
	\centering 
	\begin{tabular}{|c|c|c|c|c|c|c|c|}
		\hline  
		$f$& $(\aaa,\bbb,\ccc,\ddd)$&$f$& $(\aaa,\bbb,\ccc,\ddd)$&$f$& $(\aaa,\bbb,\ccc,\ddd)$&$f$& $(\aaa,\bbb,\ccc,\ddd)$ \\
		\hline\hline
		$8$&$(3,14,6,7)/12$&$12$&$(2,10,3,6)/9$&$12$&$(1,21,5,8)/15$&$12$&$(4,9,5,17)/15$\\
		\hline 
		$12$&$(9,28,10,23)/30$&$12$&$(3,16,10,41)/30$&$12$&$(1,2,1,3)/3$&$12$&$(3,22,6,11)/18$\\ 
		\hline
		$12$&$(3,10,6,23)/18$&$16$&$(1,3,1,4)/4$&$16$&$(1,10,2,5)/8$&$16$&$(3,14,6,31)/24$\\
		\hline
		$20$&$(5,32,6,23)/30$&$20$&$(1,16,6,43)/30$&$20$&$(1,4,1,5)/5$&$20$&$(3,38,6,19)/30$\\
		\hline
		$20$&$(1,6,2,13)/10$&$24$&$(1,5,1,6)/6$&$24$&$(3,46,6,23)/36$&$24$&$(3,22,6,47)/36$\\
		\hline
		$60$&$(1,14,1,15)/15$&$60$&$(3,118,6,59)/90$&$60$&$(3,58,6,119)/90$&&\\
		\hline
	\end{tabular}
\end{table}
\subsubsection*{Case 15.   resultant$(P(x,y,z),P(-x^2,-y^2,-z^2),z)$}

$ \left(z^{2}+1\right) \left(y -1\right) \left(y +1\right) \left(y z -1\right) \left(y^{2}+y z +z^{2}\right) \left(y^{3} z -1\right)=0$

\begin{table}[H]        
	\centering 
	\begin{tabular}{|c|c|c|c|c|c|c|c|}
		\hline  
		$f$& $(\aaa,\bbb,\ccc,\ddd)$&$f$& $(\aaa,\bbb,\ccc,\ddd)$&$f$& $(\aaa,\bbb,\ccc,\ddd)$&$f$& $(\aaa,\bbb,\ccc,\ddd)$ \\
		\hline\hline
		$6$&$(6,3,4,3)/6$&$6$&$(1,8,4,3)/6$&$6$&$(12,4,6,2)/9$&$8$&$(3,14,6,7)/12$\\
		\hline 
		&$(4,f-4,4,f)/f$&&$(6,4f-4,12,2f-2)/f$&&$(6,2f-4,12,4f-2)/f$&&\\ 
		\hline
	\end{tabular}
\end{table}

In summary, there are three infinite sequences of rational $a^3b$-quadrilaterals, and many sporadic rational solutions are merged into these three infinite sequences. The following table combines all above solutions, which all admit the simplest $2$-layer earth map tilings $T(f\aaa\bbb\ddd,2\ccc^{\frac{f}{2}})$ for some even integers $f\ge6$, together with their various modifications when $\beta$ is an integer multiple of $\gamma$ (see Table \ref{Tab-1.1} and \ref{Tab-1.2}). The first two pictures of Figure \ref{fig 1-1-1} are two simple examples with just $10$ tiles.

\begin{table}[htp]        
	\centering 
	\caption{Tilings with $\aaa\bbb\ddd$ or $\aaa\ccc\ddd$.}\label{abd}
	\begin{tabular}{|c|c|c|}
		\hline  
		Vertex& $(\aaa,\bbb,\ccc,\ddd),f$ &$(\aaa,\bbb,\ccc,\ddd),f$\\
		\hline\hline
		\multirow{8}{*}{$\aaa\bbb\ddd$}&$(6,3,4,3)/6,6$&$(1,8,4,3)/6,6$\\
		\cline{2-3}
		&$(12,4,6,2)/9,6$ &$(2,10,3,6)/9,12$\\
		\cline{2-3}
		&$(1,21,5,8)/15,12$&$(4,9,5,17)/15,12$\\  
		\cline{2-3}
		&$(9,28,10,23)/30,12$&$(3,16,10,41)/30,12$\\
		\cline{2-3}
		&$(3,20,4,13)/18,18$&$(5,32,6,23)/30,20$\\
		\cline{2-3}
		&$(1,16,6,43)/30,20$ &$(1,42,4,17)/30,30$\\
		\cline{2-3}
		&$(4,f-4,4,f)/f$&$(6,4f-4,12,2f-2)/3f$\\
		\cline{2-3}
		&$(6,2f-4,12,4f-2)/3f$&\\
		\hline  
	\end{tabular}
\end{table}

 \subsection{Tilings with two degree $3$ vertex types (excluding $\aaa\bbb\ddd$ and $\aaa\ccc\ddd$).}
 
 By  Lemma \ref{geometry2} and \ref{tildev20}, we obtain $7$ distinct cases.

      \subsubsection*{Case $\{\aaa^2\bbb,\bbb^3\}$, $(x={\mathrm e}^{\mathrm{i}\ddd},y={\mathrm e}^{\frac{8\mathrm{i} \pi}{f}})$}

   $\zeta_{3} x^{3}-\zeta_{3}^{2} x^{2} y-\left(\zeta_{3}^{2}-1\right) x^{2}+\left(\zeta_{3}-\zeta_{3}^{2}\right) x y -\zeta_{3}^{2} x +y=0$

No good rational angle solution (or simply ``None'' hereafter).

   \subsubsection*{Case $\{\aaa^2\bbb,\bbb^2\ccc\}$, $(x={\mathrm e}^\frac{\mathrm{2i} \ddd}{3},y={\mathrm e}^{\frac{\mathrm{4i} \pi}{3f}})$}

    $\zeta_{3} x^{4} y^{2}-\zeta_{3}^{2} x^{2} y^{6}+\zeta_{3} xy^{5}-2 \zeta_{3}^{2} x^{2} y^{3}+x^{3} y -\zeta_{3}^{2} x^{2}+y^{4}=0$

    None.
    
   \subsubsection*{Case $\{\aaa^2\bbb,\ccc^3\}$, $(x={\mathrm e}^{\mathrm{2i} \ddd},y={\mathrm e}^{\frac{4\mathrm{i} \pi}{f}})$}

    $-\zeta_{3}^{2} y^{4}+\left(\zeta_{3}-\zeta_{3}^{2}\right) x^{2} y -\zeta_{3}^{2} xy^{2}-\left(\zeta_{3}^{2}-1\right) y^{3}-\zeta_{3}^{2} x^{2}=0$

None.

\subsubsection*{Case $\{\aaa^2\bbb,\ccc\ddd^2\}$, $(x={\mathrm e}^{2\mathrm{i} \ddd},y={\mathrm e}^{\frac{4\mathrm{i} \pi}{f}})$}

    $xy^{4}+xy^{3}-xy^{2}-y^{3}-x y -y^{2}+y +1=0$

None.

\subsubsection*{Case $\{\bbb\ddd^2,\bbb^2\ccc\}$, $(x={\mathrm e}^\frac{\mathrm{i} \ddd}{2},y={\mathrm e}^{\frac{2\mathrm{i} \pi}{f}})$}

$x^{6}-x^{5} y -x^{5}+2 x^{3}y -x y^{2}-x y +y^{2}=0$

None.

\subsubsection*{Case $\{\bbb\ddd^2,\ccc^3\}$, $(x={\mathrm e}^\frac{\mathrm{i} \ddd}{2},y={\mathrm e}^{\frac{4\mathrm{i} \pi}{f}})$}

$-\zeta_{3}^{2} x^{2} y^{2}-\left(\zeta_{3}-1\right) x^{2} y +\zeta_{3}^{2} xy^{2}+\zeta_{3}^{2} x +\left(\zeta_{3}-1\right) y -\zeta_{3}^{2}=0$

None.

\subsubsection*{Case $\{\bbb\ddd^2,\bbb\ccc^2\}$, $(x={\mathrm e}^\frac{\mathrm{i} \ddd}{2},y={\mathrm e}^{\frac{4\mathrm{i} \pi}{f}})$}

$x^{5} y -x^{4} y +x^{3} y^{2}-x^{4}-xy^{2}+x^{2}-x y +y=0$

$f=16,(\aaa,\bbb,\ccc,\ddd)=(1,4,2,2)/4$.

 \vspace{0.5cm}
 
 The above results are summarized in Table \ref{3ab1}, and there is only one rational $a^3b$-monotile admitting two tilings with $16$ tiles, shown in Table \ref{Tab-1.1} and the third, fourth pictures of Figure \ref{fig 1-1-1}. 
 \begin{table}[htp]        
 	\centering 
 	\caption{Tilings with two degree $3$ vertex types (excluding $\aaa\bbb\ddd$ and $\aaa\ccc\ddd$).}
 	\label{3ab1}
 	\begin{tabular}{|c|c|c|}
 		\hline  
 		Vertex & $(\aaa,\bbb,\ccc,\ddd),f$ & Contradiction\\
 		\hline\hline
 		$\aaa^2\bbb,\bbb^3$& & \multirow{7}{*}{ No solution }\\
 		\cline{1-2}
 		$\aaa^2\bbb,\bbb^2\ccc$& & \\
 		\cline{1-2}
 		$\aaa^2\bbb,\ccc^3$& & \\
 		\cline{1-2}
 		$\aaa^2\bbb,\ccc\ddd^2$& & \\
 		\cline{1-2}
 		$\bbb\ddd^2,\bbb^2\ccc$& & \\
 		\cline{1-2}
 		$\bbb\ddd^2,\ccc^3$& & \\
 		\cline{1-3}
 		$\bbb\ddd^2,\bbb\ccc^2$& $(1,4,2,2)/4,16$& \\
 		\hline
 	\end{tabular}      	
 \end{table}

 \subsection{Tilings with a unique degree $3$ vertex type (excluding $\aaa\bbb\ddd$ and $\aaa\ccc\ddd$)}

By  Lemma \ref{combined_lemma}, we obtain $5+5+9+9=28$ distinct cases.

\subsubsection*{Case $\{\bbb\ccc^2,\aaa^4\}$, $(x={\mathrm e}^{\mathrm{i} \ccc},y={\mathrm e}^{\frac{4\mathrm{i} \pi}{f}})$}

$x^{5} y^{2}-x^{4} y^{2}-x^{4} y -x^{3} y +x^{2}y+x y +x -1=0$

None.

\subsubsection*{Case $\{\bbb\ccc^2,\aaa^3\ddd\}$, $(x={\mathrm e}^\frac{2\mathrm{i} \ddd}{3},y={\mathrm e}^{\frac{4\mathrm{i} \pi}{f}})$}

$\zeta_{3} x^{5} y +x^{5} \zeta_{3}-x^{6} \zeta_{3}^{2}-x^{3} y^{2} \zeta_{3}^{2}-x^{3} y \zeta_{3}^{2}-y^{3} \zeta_{3}^{2}+x y^{3}+xy^{2}=0$

$f=36,(\aaa,\bbb,\ccc,\ddd)=(5,4,7,3)/9$.

\subsubsection*{Case $\{\bbb\ccc^2,\aaa^2\ddd^2\}$, $(x={\mathrm e}^{2\mathrm{i} \ddd},y={\mathrm e}^{\frac{4\mathrm{i} \pi}{f}})$}

$x y^{2}+2 y^{3}-x y +y^{2}-2 x -y=0$

None.
\subsubsection*{Case $\{\bbb\ccc^2,\aaa\ddd^3\}$, $(x={\mathrm e}^{\mathrm{2i} \ddd},y={\mathrm e}^{\frac{4\mathrm{i} \pi}{f}})$}

$x^{4} y^{2}+x^{3} y^{3}-x^{2} y^{3}-x^{2} y^{2}-x^{2} y -x^{2}+x +y=0$

$f=10,(\aaa,\bbb,\ccc,\ddd)=(1,6,2,3)/5,10$.

\subsubsection*{Case $\{\bbb\ccc^2,\ddd^4\}$, $(x={\mathrm e}^{\mathrm{i} \ccc},y={\mathrm e}^{\frac{4\mathrm{i} \pi}{f}})$}

$x^{3} y^{2}-x^{3} y +x^{2} y -xy^{2}+x^{2}-x y +y -1=0$

$f=16,(\aaa,\bbb,\ccc,\ddd)=(1,4,2,2)/4$.

  \subsubsection*{Case $\{\bbb^3,\aaa^4\}$, $(x={\mathrm e}^{2\mathrm{i} \ddd},y={\mathrm e}^{\frac{8\mathrm{i} \pi}{f}})$}

$\zeta_{12}^{4} x^{3}+\zeta_{12}^{3} x^{2} y -\left(\zeta_{12}^{5}-\zeta_{12}\right) x^{2}+\left(\zeta_{12}^{4}-1\right) x y +\zeta_{12}^{2} x +\zeta_{12} y=0$

$f=6,(\aaa,\bbb,\ccc,\ddd)=(3,4,8,1)/6$.

$f=8,(\aaa,\bbb,\ccc,\ddd)=(3,4,6,2)/6$.

   \subsubsection*{Case $\{\bbb^3,\aaa^3\ddd\}$, $(x={\mathrm e}^\frac{\mathrm{2i} \ddd}{3},y={\mathrm e}^{\frac{4\mathrm{i} \pi}{f}})$}

$\zeta_{3} x^{4}-\zeta_{3}^{2} x^{3} y +\left(\zeta_{3}-\zeta_{3}^{2}\right) x^{2} y -\left(\zeta_{3}^{2}-1\right) x^{2}-\zeta_{3}^{2} x +y=0$

None.
  \subsubsection*{Case $\{\bbb^3,\aaa^2\ddd^2\}$, $(x={\mathrm e}^{2\mathrm{i} \ddd},y={\mathrm e}^{\frac{4\mathrm{i} \pi}{f}})$}

$-\left(2 \zeta_{3}^{2}-1\right) x y -x + y -1+2 \zeta_{3}=0$

None.

  \subsubsection*{Case $\{\bbb^3,\aaa\ddd^3\}$, $(x={\mathrm e}^{\mathrm{2i} \ddd},y={\mathrm e}^{\frac{4\mathrm{i} \pi}{f}})$}

$\left(\zeta_{3}-\zeta_{3}^{2}\right) x^{4} y -\zeta_{3} x^{3} y +\zeta_{3}^{2} x^{2} y +\zeta_{3}^{2} x^{2}-x -\zeta_{3}^{2}+1=0$

None.

  \subsubsection*{Case $\{\bbb^3,\ddd^4\}$, $(x={\mathrm e}^{\mathrm{i} \ccc},y={\mathrm e}^{\frac{4\mathrm{i} \pi}{f}})$}

$-\left(\zeta_{3}^{2}-1\right) x^{3}+\left(\zeta_{3}+\zeta_{3}^{2}\right) x^{2} y -\left(\zeta_{3}+\zeta_{3}^{2}\right) x y -\left(1-\zeta_{3}\right) y^{2}=0$

None.

   \subsubsection*{Case $\{\aaa^2\bbb,\bbb^4\}$, $(x={\mathrm e}^{\mathrm{i} \ddd},y={\mathrm e}^{\frac{4\mathrm{i} \pi}{f}})$}

$\zeta_{8}^{3} x^{3}+\zeta_{8}^{2} x^{2} y +\left(\zeta_{8}^{2}+1\right) x^{2}+\left(\zeta_{8}^{3}+\zeta_{8}\right) x y +\zeta_{8} x +y=0$

$f=6,(\aaa,\bbb,\ccc,\ddd)=(9,6,12,5)/12$.

   \subsubsection*{Case $\{\aaa^2\bbb,\bbb^3\ccc\}$, $(x={\mathrm e}^\frac{2\mathrm{i} \ddd}{5},y={\mathrm e}^{\frac{4\mathrm{i} \pi}{5f}})$}

$\zeta_{5}^{2} x^{7} y^{2}-\zeta_{5}^{4} x^{4} y^{8}+\zeta_{5}^{2} x^{2} y^{7}-\zeta_{5}^{4} x^{4} y^{3}-\zeta_{5}^{3} x^{3} y^{5}+x^{5} y -\zeta_{5}^{3} x^{3}+y^{6}=0$

None.

\subsubsection*{Case $\{\aaa^2\bbb,\bbb^2\ccc^2\}$, $(x={\mathrm e}^{2\mathrm{i} \ddd},y={\mathrm e}^{\frac{4\mathrm{i} \pi}{f}})$}

$y^{4}+xy^{2}+2 y^{3}-2 x y -y^{2}-x=0$

None.

\subsubsection*{Case $\{\aaa^2\bbb,\bbb^2\ddd^2\}$, $((x={\mathrm e}^\frac{\mathrm{i} \ddd}{2},y={\mathrm e}^{\frac{4\mathrm{i} \pi}{f}})$}

$\zeta_{4} x^{7}+x^{6} y +x^{6}+x^{4}+\zeta_{4} x^{3} y +\zeta_{4} x y +\zeta_{4} x +y=0$

None.

\subsubsection*{Case $\{\aaa^2\bbb,\bbb\ccc^3\}$, $(x={\mathrm e}^{2\mathrm{i} \ddd},y={\mathrm e}^{\frac{4\mathrm{i} \pi}{f}})$}

$-x^{2} y^{8}+y^{10}-x^{3} y^{6}-x^{2} y^{7}+xy^{8}+y^{9}+x^{6} y +x^{5} y^{2}-x^{4} y^{3}-x^{3} y^{4}+x^{6}-x^{4} y^{2}=0$

$f=36,(\aaa,\bbb,\ccc,\ddd)=(15,6,10,7)/18$.

\subsubsection*{Case $\{\aaa^2\bbb,\bbb\ccc\ddd^2\}$, $(x={\mathrm e}^{\mathrm{2i} \ddd},y={\mathrm e}^{\frac{4\mathrm{i} \pi}{f}})$}

$x^{2} y^{4}+x^{2} y^{3}-x^{2} y^{2}-x y^{3}-x y -y^{2}+y +1=0$

None.

\subsubsection*{Case $\{\aaa^2\bbb,\ccc^4\}$, $(x={\mathrm e}^{\mathrm{2i} \ddd},y={\mathrm e}^{\frac{4\mathrm{i} \pi}{f}})$}

$\left(\zeta_{4}^{2}+\zeta_{4}\right) y^{4}+2 \zeta_{4} x^{2} y +2 \zeta_{4} x y^{2}+2 \zeta_{4} y^{3}+\left(\zeta_{4}+1\right) x^{2}=0$

None.

\subsubsection*{Case $\{\aaa^2\bbb,\ccc^2\ddd^2\}$, $(x={\mathrm e}^{\mathrm{i} \ddd},y={\mathrm e}^{\frac{2\mathrm{i} \pi}{f}})$}

$-xy^{4}+x^{3} y -x^{2} y^{2}-x y^{3}+x^{2} y +x y^{2}-y^{3}+x^{2}
=0$

None.
\subsubsection*{Case $\{\aaa^2\bbb,\ddd^4\}$, $(x={\mathrm e}^{\mathrm{i} \ccc},y={\mathrm e}^{\frac{2\mathrm{i} \pi}{f}})$}

$x^{3} y +x^{2} y^{2}-x y^{3}+y^{4}+x^{3}-x^{2} y +x y^{2}+y^{3}=0$

None.

  \subsubsection*{Case $\{\bbb\ddd^2,\aaa^4\}$, $(x={\mathrm e}^{\mathrm{i} \ccc},y={\mathrm e}^{\frac{4\mathrm{i} \pi}{f}})$}

$x^{5}+x^{4} y +x^{2} y^{3}-xy^{4}-x^{4}+x^{3} y +xy^{3}+y^{4}=0$

None.

  \subsubsection*{Case $\{\bbb\ddd^2,\aaa^2\bbb^2\}$, $(x={\mathrm e}^{\mathrm{i} \ccc},y={\mathrm e}^{\frac{4\mathrm{i} \pi}{f}})$}

$x^{7}-x^{6} y -x^{5} y^{2}+x^{2} y^{5}+x^{5} y -x^{2} y^{4}-xy^{5}+y^{6}=0$

None.

  \subsubsection*{Case $\{\bbb\ddd^2,\aaa^2\bbb\ccc\}$, $(x={\mathrm e}^{\mathrm{2i} \ddd},y={\mathrm e}^{\frac{4\mathrm{i} \pi}{f}})$}

$ x^{2}y^{3}-xy^{4}  +y^{4}-x^{2}y -y^{3}+x^{2}-x +y=0$

$f=16,(\aaa,\bbb,\ccc,\ddd)=(1,4,2,2)/4$.

  \subsubsection*{Case $\{\bbb\ddd^2,\aaa^2\ccc^2\}$, $(x={\mathrm e}^{\mathrm{i} \ccc},y={\mathrm e}^{\frac{4\mathrm{i} \pi}{f}})$}

$\left(y -1\right) \left(y +1\right) \left(x -1\right) \left(x +1\right)\! \left(x^{2} y -xy^{2}+x y -x +y \right)=0$

$f=8,(\aaa,\bbb,\ccc,\ddd)=(2,6,4,3)/6$.

  \subsubsection*{Case $\{\bbb\ddd^2,\bbb^4\}$, $(x={\mathrm e}^{\mathrm{i} \ccc},y={\mathrm e}^{\frac{4\mathrm{i} \pi}{f}})$}

$\left(\zeta_{4}+1\right) x^{3}-2 x^{2} y +2 x y +\left(\zeta_{4}-1\right) y^{2}=0$

$f=8,(\aaa,\bbb,\ccc,\ddd)=(11,6,4,9)/12$.

  \subsubsection*{Case $\{\bbb\ddd^2,\bbb^3\ccc\}$, $(x={\mathrm e}^\frac{\mathrm{i} \ccc}{3},y={\mathrm e}^{\frac{4\mathrm{i} \pi}{f}})$}

$\zeta_{3} x^{7} y -\zeta_{3}^{2} x^{9}-\zeta_{3} x^{4} y +x^{8}+\zeta_{3} xy^{2}-x^{5} y -\zeta_{3}^{2} y^{2}+x^{2} y=0$

None.

  \subsubsection*{Case $\{\bbb\ddd^2,\bbb^2\ccc^2\}$, $(x={\mathrm e}^{\mathrm{i} \ccc},y={\mathrm e}^{\frac{4\mathrm{i} \pi}{f}})$}

$x^{3}y+x^{3}-x^{2} y +xy^{2}  +x^{2}-x y +y^{2}+y=0$

None.

  \subsubsection*{Case $\{\bbb\ddd^2,\bbb\ccc^3\}$, $(x={\mathrm e}^\frac{2\mathrm{i} \ddd}{3},y={\mathrm e}^{\frac{4\mathrm{i} \pi}{f}})$}

$(x - 1)(x + 1)(x^2 + x + 1)(x^2 - x + 1)(x^4y - x^3y + x^2y^2 - x^2 + xy - y)
=0$

None.

 \subsubsection*{Case $\{\bbb\ddd^2,\ccc^4\}$, $(x={\mathrm e}^{\mathrm{2i} \ddd},y={\mathrm e}^{\frac{4\mathrm{i} \pi}{f}})$}

$\zeta_{4} x^{2} y^{2}-\left(\zeta_{4}-1\right) x^{2} y -\zeta_{4} xy^{2}-x +\left(\zeta_{4}-1\right) y +1=0$

$f=16,(\aaa,\bbb,\ccc,\ddd)=(1,4,2,2)/4$.

\vspace{0.5cm}

 The above results are summarized in Table \ref{3ab2}. There are only two new rational $a^3b$-monotiles each admitting a unique tiling with $36$ tiles, shown in Table \ref{Tab-1.1}  and the last two pictures of Figure \ref{fig 1-1-1}.

\begin{table}[htp]        
	\centering 
	\caption{Tilings with a unique degree $3$ vertex type (excluding $\aaa\bbb\ddd$ and $\aaa\ccc\ddd$)}\label{3ab2}
	\begin{tabular}{|c|c|c|c|}
		\hline  
		\multicolumn{2}{|c|}{Vertex}& $(\aaa,\bbb,\ccc,\ddd),f$ & Contradiction\\
		\hline\hline
		\multirow{4}{*}{$\bbb\ccc^2$} & $\aaa^4,\aaa^2\ddd^2$ & & No solution \\
		\cline{2-4}
		&$\aaa^3\ddd$&$(5,4,7,3)/9,36$ &\\
		\cline{2-4}
		&$\aaa\ddd^3$&$(1,6,2,3)/5,10$& already in Table \ref{abd}\\
		\cline{2-4}
		&$\ddd^4$& $(1,4,2,2)/4,16$& already in Table \ref{3ab1}\\
		\hline
		\hline
		\multirow{3}{*}{$\bbb^3$} 
		&\multirow{2}{*}{$\aaa^4$}& $(3,4,8,1)/6,6$
		&already in Table \ref{abd} ($\aaa\leftrightarrow\ddd$ and $\bbb\leftrightarrow\ccc$)\\
		\cline{3-4}
		&&$(3,4,6,2)/6,8$&no tiling by Balance Lemma\\
		\cline{2-4}		
		&$\aaa^3\ddd,\aaa^2\ddd^2,\aaa\ddd^3,\ddd^4$&& No solution\\		
		\hline
		\hline				
		\multirow{3}{*}{$\aaa^2\bbb$} 			
		&$\bbb^4$& $(9,6,12,5)/12,6$&no tiling by Balance Lemma\\
		\cline{2-4}	
		&$\bbb\ccc^3$&$(15,6,10,7)/18,36$&\\	
		\cline{2-4}		
		&$\bbb^3\ccc,\bbb^2\ccc^2,\bbb^2\ddd^2,\bbb\ccc\ddd^2,\ccc^4,\ccc^2\ddd^2,\ddd^4$&&No solution\\
		\hline
		\hline
		\multirow{3}{*}{$\bbb\ddd^2$} & $\aaa^4,\aaa^2\bbb^2,\bbb^4,\bbb^3\ccc,\bbb^2\ccc^2,\bbb\ccc^3$&&No solution \\
		\cline{2-4}		
		&$\aaa^2\bbb\ccc,\ccc^4$&$(1,4,2,2)/4,16$&already in Table \ref{3ab1}\\
		\cline{2-4}	
		&$\aaa^2\ccc^2$&$(2,6,4,3)/6,8$&no tiling by Balance Lemma\\
		\hline
	\end{tabular} 
\end{table}

\end{document}